\documentclass[12pt,oneside,a4paper,leqno]{article}

\usepackage{amsfonts,amssymb}
\usepackage[centertags]{amsmath}
\usepackage[all]{xy}\xyoption{all}
\usepackage{amsthm}
\newcounter{lemma}
\swapnumbers
\theoremstyle{plain}
\newtheorem{lemma}[equation]{Lemma}\numberwithin{lemma}{section}
\newtheorem{theo}[equation]{Theorem}

\newtheorem{coro}[equation]{Corollary}
\theoremstyle{definition}
\newtheorem{defi}[equation]{Definition}
\newtheorem{example}[equation]{Example}
\newtheorem{remark}[equation]{Remark}
\theoremstyle{remark}

\numberwithin{equation}{section}

\makeatletter
\def\tagform@#1{\maketag@@@{\ignorespaces#1\unskip\@@italiccorr}}

\makeatother

\newcommand{\ldiag}[2]{%
\begin{equation}\label{#1}\begin{aligned}\xymatrix{#2}\end{aligned}\end{equation}%
}
\newcommand{\diag}[1]{%
\begin{equation}\begin{aligned}\xymatrix{#1}\end{aligned}\end{equation}%
}
\newcommand{\ndiag}[1]{%
\begin{equation*}\begin{aligned}\xymatrix{#1}\end{aligned}\end{equation*}%
}

\newdir{ >}{!/8pt/@{ }*@{>}}
\newdir{< }{!/-8pt/@{ }*@{<}}

\theoremstyle{plain}
\newtheorem{pbtheo}[equation]{Pull back theorem}
\newtheorem{bundletheo}[equation]{Bundle theorem}

\newcommand{\hh}{\Phi}

\newcommand{\vu}{\nu}
\newcommand{\minus}{\smallsetminus\/}

\newcommand{\R}{\mathbb{R}}
\newcommand{\vect}{{\mathbf{Vect}}}
\renewcommand{\vec}{{\mathbf{Vect}}}
\renewcommand{\Vec}{{\mathbf{Vect}}}

\newcommand{\V}{{\mathcal{{\bf V}}}}

\newcommand{\W}{{\mathcal{{\bf W}}}}
\newcommand{\K}{\mathbb{K}}
\newcommand{\KK}{\mathcal{K}}
\newcommand{\Ab}{\mathbf{Ab}}

\newcommand{\Vtop}{\mbox{$\V$-${{{\mathbf{Top}}}}$}}

\newcommand{\homotopic}{\sim}
\renewcommand{\top}{\mathbf{Top}} 

\newcommand{\T}{\mathbb{T}}
\newcommand{\point}{\circ}
\newcommand{\VCW}{\mbox{$\V$-CW}}
\newcommand{\bVCW}{\mbox{$\bar \V$-CW}}
\newcommand{\WCW}{\mbox{$\W$-CW}}
\newcommand{\bWCW}{\mbox{$\bar \W$-CW}}

\newcommand{\op}{\mathrm{op}}

\newcommand{\from}{\colon}
\newcommand{\st}{\ \mathrm{:} \ }
\newcommand{\VpCW}{\mbox{$\V'$-CW}}

\newcommand{\Stra}{\mathbf{Stra}}
\newcommand{\stra}{\mathbf{Stra}}
\newcommand{\VStra}{\mbox{$\V$-$\mathbf{Stra}$}}
\newcommand{\WStra}{\mbox{$\W$-$\mathbf{Stra}$}}
\newcommand{\Vstra}{\mbox{$\V$-$\mathbf{Stra}$}}
\newcommand{\Vpstra}{\mbox{$\V'$-$\mathbf{Stra}$}}
\newcommand{\Wstra}{\mbox{$\W$-$\mathbf{Stra}$}}
\newcommand{\vvect}{\Vstra}
\newcommand{\Vvect}{\Vstra}

\newcommand{\ff}{F}

\begin{document}
\pagenumbering{arabic}

\title{%
$K$-theory of stratified vector bundles
}

\author{Hans-Joachim Baues and Davide L.~Ferrario}



\date{}
\maketitle

\begin{abstract}
We show that the Atiyah-Hirzebruch $K$-theory of spaces admits 
a canonical generalization for stratified spaces.
For this we study algebraic constructions on stratified
vector bundles. In particular the tangent bundle
of a stratified manifold is such a stratified 
vector bundle.
\end{abstract}




\section{Families of vector spaces}
Let  $\R$ be the field of real numbers 
and $\R^n = \R \oplus \dots \oplus \R$ be the standard 
$n$-dimensional $\R$-vector space.
Let $\vect$ denote the category of finite dimensional 
$\R$-vector spaces  
and linear maps and  $\V$ a subcategory of $\vect$,
termed the \emph{structure category}.
For example let $\V$ be the subcategory of the surjective
maps in $\Vec$. Or let $G$ be a 
subgroup of the automorphism
group $GL_n(\R)$ of $\R^n$. Then $G$ yields the subcategory 
$G\subset \vec$ consisting of one object $\R^n$ and morphisms
given by elements in $G$.

The category $\vect$ is a topological category. We say 
that $\V$ is a \emph{closed} subcategory
if for objects $V$, $W\in \V$ the inclusion
of morphisms sets
\begin{equation*}
\hom_\V(V,W) \subset \hom_{\vect}(V,W)
\end{equation*}
is a closed subspace.
Here $\hom_\vect(V,W)$ is homeomorphic to $\R^N$ 
while $\hom_{\V}(V,W)$ 
needs not to be a vector space.

Following Atiyah \cite{at} a \emph{family of vector spaces} 
in $\V$
is a topological space $X$ 
together with a map $p_X:X \to \bar X$ and for every $b\in \bar X$ 
a homeomorphism, named \emph{chart},
$\Phi_b: p^{-1}_X(b)  \approx X_b $ where $X_b$ is 
a suitable object of $\V$ depending upon $b\in \bar X$. 
We denote the family $(p_X:X \to \bar X, X_b, \Phi_b, b\in \bar X)$
simply by $X$ (it will be clear from the context whether $X$ 
is the underlying space or the family).
The map $p_X$ is termed the \emph{projection}, the space $X$ the 
\emph{total space} of the family of vector spaces, 
the space $\bar X$ the \emph{base space}, and for every $b$ the vector
space $X_b$ is termed \emph{fiber} over $b$. A family $X$ is also
termed $\V$-family, to make explicit the choice of $\V$.
The $\V$-family $X$ is \emph{discrete} if $\bar X$ is a discrete space.
In this case $X$ is a disjoint union of objects in $\V$
and we call $X$ a \emph{$\V$-set}.

Given two $\V$-families $X$
and $Y$ a \emph{$\V$-map} from $X$ to $Y$ 
(equivalently, a \emph{homomorphism} or an \emph{admissible
map controlled by $\V$})
is a pair of maps
$(f,\bar f)$ 
such that the following diagram commutes
\diag{%
X \ar[r]^f
\ar[d]_{p_X} & Y \ar[d]^{p_Y} \\
\bar X \ar[r]^{\bar f} & \bar Y,
}%
and such that for every $b\in \bar X$ the composition
given by the dashed arrow of the diagram
\ndiag{%
p_X^{-1}(b)\ar[r]^{f|p^{-1}(b)}
&  p_Y^{-1}(\bar f b) \ar[d]^{\Phi_{\bar f b}} \\
X_b 
\ar[u]^{\Phi_b^{-1}} 
\ar@{-->}[r]^{} & Y_{\bar f b}
}%
is a morphism of $\V$.
A $\V$-isomorphism $(f,\bar f)\from X \to Y$
induces homeomorphisms $f\from X\to Y$, $\bar f\from \bar X \to \bar Y$,
and $f|X_b\from X_b  \to Y_{\bar f b}$ for every $b\in \bar X$.

If $V$ is a vector space in $\V$ and $\bar X$ is a topological space,
then 
the projection onto the first factor
$p_1\colon X = \bar X \times V \to \bar X$
is the projection map of the \emph{product family} with fiber $V$;
the charts $\Phi_b: \{b\}\times V \to V=X_b$ are given by projection 
 and $X_b=V$ for all $b\in \bar X$.

If $Y$ is a family  $\V$-isomorphic  to a product family then
$Y$ is said to be a \emph{trivial} vector bundle.

Given a family $Y$  and 
a map $\bar f\colon \bar X \to \bar Y$, the pull-back $X=\bar f^*Y$
is 
defined by the following pull-back diagram.
\diag{%
X =
\bar f^* Y 
\ar@{-->}[r]
\ar@{-->}[d]
&  Y
\ar[d]^{p}
\\
\bar X 
\ar[r]^{\bar f}
& \bar Y.\\
}%
If $\bar f = i\colon \bar X \subset \bar Y$ is an inclusion, then
$i^*Y$ is called the \emph{restriction} of $Y$ to $\bar X$
and is denoted by $Y|\bar X$. 

A \emph{$\V$-vector bundle} or a $\V$-bundle
is a locally trivial family of vector 
spaces in $\V$, i.e.  a family $X$ over $\bar X$ such that 
every $b\in \bar X$ admits a neighborhood $U$ for which $X|U$ is trivial.

If $X$ is a $\V$-family  and $Z$ is a 
topological space, then the cartesian product 
$X\times Z$ is a $\V$-family
with projection 
$p_{X\times Z} = p\times 1_Z\colon X\times Z \to \bar X \times Z$.
The fiber over  a point $(b,z)\in \bar X\times Z$ is equal to 
$p_{X}^{-1}(b)\times \{z\}$ and 
the chart $\Phi_{(b,z)}$ is defined by $(x,z) \mapsto \Phi_b(x) \in X_b$,
where of course we set $(X\times Z)_{(b,z)} = X_b$ and $\Phi_b$ 
is 
the chart $\Phi_b:p_{X}^{-1}(b) \to X_b$.

In particular, if $Z$ is  the unit interval $I$ 
one obtains the cylinder object $X\times I$;
this leads to the definition of \emph{homotopy}: 
two $\V$-maps
$f_0,f_1\colon X \to Y$ 
are $\V$-homotopic (in symbols $f_0\homotopic f_1$)
if there is a $\V$-map $F\colon X\times I \to Y$ 
such that $f_0=Fi_0$ and $f_1=Fi_1$. Here $i_0$ and $i_1$ 
are the inclusions $X \to X\times I$ 
to $X\times \{0\}$ and $X\times \{1\}$ respectively.

Let $\top$ be the category of 
topological spaces 
and (continuous) maps. Then $\Vtop$ is the category 
of $\V$-families $p_X: X \to \bar X$ in $\top$ and $\V$-maps.
Homotopy of $\V$-maps yields 
the homotopy relation $\sim$ which is 
a natural equivalence relation 
on $\Vtop$ so that the homotopy category
$(\Vtop)/_\homotopic$ is defined.

A $\V$-map $i:A \to Y$ is termed a \emph{closed
inclusion} if $\bar i: \bar A \to \bar Y$ is an inclusion, 
$\bar i \bar A$ is closed in $\bar Y$ and  $\bar i^*Y = A$.
A closed inclusion
$i: A \to Y$ induces isomorphisms 
on fibers.

The push-out construction can be extended to the category $\Vtop$,
provided the push-out is defined by a closed inclusion. 
\begin{lemma}\label{lemma:pushout}
Given $\V$-families $A$, $X$, $Y$ and $\V$-maps $f\from A \to X$,
$i\from A\to Y$ with $i$ a closed inclusion,
the push-out diagram 
in $\Vtop$ 
\diag{%
A \ar[r]^f \ar@{->}[d]_{i} & X \ar@{->}[d]\\
Y \ar[r] & Z \\
}%
exists 
and 
$X\to Z$ is a closed inclusion. 
\end{lemma}
\begin{proof}
See \cite{qq}, lemma (2.5). 
\end{proof}

\section{Stratified vector bundles}

We introduce the notion of stratified space and stratified vector 
bundle. 
They are particular cases of stratified bundles as defined in
\cite{qq}, definition (4.1). 
In the next section we describe as an example the stratified
tangent bundle.

Let $A$ be a closed subset  of a space $X$. 
We say that $(X,A)$ is a \emph{CW-pair} if there exists a homeomorphism
$(X,A) \approx (X',A')$ of pairs where $X'$ is a CW-complex
and $A'$ a subcomplex of $X'$.
A \emph{CW-space} $X$ is a space homeomorphic to a CW-complex.

\begin{defi}\label{defi:strati}
A  \emph{stratified} space is a space $\bar X$ together with 
 a filtration 
\begin{equation*}
\bar X_0 \subset \bar X_1 \subset \dots
\subset \bar X_n \subset \dots \subset \lim_{n\to \infty} \bar X_n = \bar X
\end{equation*}
such that for every $i\geq 1$ there is  a
CW-pair $(\bar M_i, \bar A_i)$ 
and an attaching map
$\bar h_i\from \bar A_i \to \bar X_{i-1}$ with the property that 
the subspace $\bar X_i$ is obtained by attaching $\bar M_i$ to $\bar X_{i-1}$ 
via the attaching map $\bar h_i$, i.e.
$\bar X_i = \bar M_i \cup_{\bar h_i} \bar X_{i-1}$.
\end{defi}

The complements $\bar X_i\minus \bar X_{i-1}$ are termed \emph{strata}, 
 while the filtration $\bar X_0 \subset
\bar X_1 \subset \dots \subset \bar X_n\subset \dots$ is called \emph{stratification} of $\bar X$.
In definition \ref{defi:strati} the strata $\bar X_i\minus \bar X_{i-1}$
coincide with the complements $\bar M_i\minus \bar A_i$. 
The pairs $(\bar M_i, \bar A_i)$ are the \emph{attached spaces} of $\bar X$.
A stratified space $\bar X$ 
is Hausdorff and regular 
and if all the attaching maps $\bar h_i$ are cellular then $\bar X$ 
is a CW-complex. 
If the number of non-empty strata of $\bar X$ is finite,
then $\bar X$ is said \emph{finite}. 
We always assume that a stratified space which is not finite
has cellular attaching maps.
A finite stratified space  $\bar X$ is \emph{compactly stratified}
if all the attached spaces $\bar M_i$ and $\bar X_0$ are compact.

A map
$\bar f\from \bar X \to \bar X'$
between stratified spaces is a filtration preserving
map $\bar f=\{\bar f_n\}_{n\ge 0}$ 
together with commutative diagrams
\ndiag{%
\bar X_{i-1}
\ar[d]^{\bar f_{i-1}}
&
\ar[l]
\bar A_i
\ar@{ >->}[r]
\ar[d]
& 
\bar M_i
\ar[d]^{\bar g_i}
\\
\bar X'_{i-1}
&
\bar A'_i
\ar[l]
\ar@{ >->}[r]
&
\bar M_i'
\\
}%
such that $\bar g_i\cup \bar f_{i-1} = \bar f_i$ for $i\geq 1$.
A map is termed \emph{stratum-preserving}
if for every $i$  
\[
\bar f_i (\bar X_i \minus \bar X_{i-1} )
\subset \bar X'_i \minus \bar X'_{i-1}.
\]
Let $\Stra$
be the category of 
finite stratified spaces 
and stratum-preserving maps.
Two maps $f_0$ and $f_1$ are homotopic in $\Stra$ if there exists
a homotopy $f_t$ such that the map
$f_t$ is a map in $\Stra$ for every $t\in I$.
With this notion of homotopy 
the homotopy category 
$\Stra/_{\sim}$  is defined.

\begin{defi}\label{defi:stratifiedvectorbundle}
A \emph{$\V$-stratified bundle} 
is a stratified space $\bar X$ together with a $\V$-family
\begin{equation*}
X \to \bar X
\end{equation*}
with the following properties. For $i\geq 1$ the restriction 
$X_i = X|\bar X_i$ is the push-out of $\V$-maps
\ndiag{%
M_i
\ar[d]
&
\ar@{ >->}[l]
A_i
\ar[d]
\ar[r]^{h_i}
&
X_{i-1} 
\ar[d]
\\
\bar M_i
&
\bar A_i
\ar@{ >->}[l]
\ar[r]^{\bar h_i}
&
\bar X_{i-1} 
\\
}%
where $M_i \to \bar M_i$ is a $\V$-vector bundle and $A_i=M_i|\bar A_i$.
Moreover, $X_0 \to \bar X_0$ is a $\V$-vector bundle
and $X=\lim_{i\to \infty} X_i$.
If $(\bar M_i, \bar A_i)$ is a disjoint union of $i$-dimensional
cells $(D^i,S^{i-1})$ then $\bar X$ is a CW-complex
and $X$ is termed a \emph{$\V$-complex}.
\end{defi}

A \emph{$\V$-stratified map} $f\from X \to X'$ 
between $\V$-stratified bundles is
given by sequences
$\{f_n\}_{n\geq 0}$ and $\{g_n\}_{n\geq 1}$ of $\V$-maps
such that given the commutative diagrams 
\ndiag{%
X_{i-1}
\ar[d]^{f_{i-1}}
&
\ar[l]
A_i
\ar@{ >->}[r]
\ar[d]
& 
M_i
\ar[d]^{g_i}
\\
X'_{i-1}
&
A'_i
\ar[l]
\ar@{ >->}[r]
&
M_i'
\\
}%
for every $i$ we have 
$g_i\cup f_{i-1} = f_i$ for $i\geq 1$.
The $\V$-map $f$ is termed \emph{stratum-preserving}
if 
for every $i$ 
\[
f(X_i\minus X_{i-1}) \subset X'_i \minus X_{i-1}'.
\]
Let $\vvect$ denote the category of $\V$-stratified 
bundles and stratum-preserving $\V$-stratified maps. 
Of course a $\V$-isomorphism
in $\vvect$ is a stratum-preserving $\V$-stratified map $f$ with an
inverse $f'$.

The following theorems  are proved in \cite{qq} (theorem (4.7)
and theorem 
(5.1)).

\begin{pbtheo}\label{theo:pullback}
Let $\bar X$ and $\bar X'$ be finite stratified spaces with attached 
spaces which are locally finite and countable CW-complexes.
Let $\bar f\from \bar X \to \bar X'$
be a stratum-preserving map (in $\Stra$) and $X'\to\bar X'$
a $\V$-stratified bundle. 
Then the pull-back $\bar f^*X' \to \bar X$ is a $\V$-stratified bundle.
\end{pbtheo}

\begin{bundletheo}
\label{theo:bundle}
Let $\V\subset \vect$ be a subcategory such that all morphisms
in $\V$ are isomorphisms (i.e. $\V$ is a groupoid).
Then a $\V$-stratified bundle $X \to \bar X$ 
is a $\V$-bundle over $\bar X$. Conversely,
given a stratified space $\bar X$ and a $\V$-bundle
$X$ over $\bar X$ then $X$ is a $\V$-stratified bundle.
\end{bundletheo}
In particular, a $\V$-complex $X\to\bar X$ 
is a $\V$-bundle over $\bar X$ and,
given a CW-complex $\bar X$ and a $\V$-bundle
$X$ over $\bar X$, then $X$ is $\V$-isomorphic to a  $\V$-complex.

\section{Stratified manifolds}
In this section we describe special families of vector spaces which arise as
``stratified tangent bundles``
of stratified manifolds or of manifolds
with singularities.
Such stratified tangent bundles are motivating examples of 
stratified vector bundles above.
Concerning the general theory of stratification of manifolds we refer the 
reader for example to \cite{kreck,thom,weinberger,pf,verona}. 
We here introduce stratified manifolds only since they show
that 
stratified vector bundles
are natural generalizations of vector bundles in a similar
way as manifolds with singularities 
or stratifolds \cite{kreck}
generalize manifolds.

If $M$ is a manifold with boundary $\partial M$,
$N$ a closed manifold, and 
$h\colon \partial M \to N$
a submersion, then
the push-out of $h$ and the inclusion $\partial M \subset M$
is a stratified space $X=M \cup_h N$ with stratification
\begin{equation*}
X_0=N \subset X_1=X.
\end{equation*}

This is a stratifold in the sense of \cite{kreck}
and generalizes ``manifolds with singularities''.
See Rudyak \cite{rudyak}, Baas \cite{baas},
Botvinnik \cite{botvinnik},
Sullivan \cite{sullivan},
Vershinin \cite{vershinin}.
Below we introduce differentiable stratified spaces which we consider as 
the general form of stratified manifolds

If $M$ is a differentiable manifold then the
tangent bundle $TM\to M$ restricted to the boundary
$TM_{|\partial M}$ is isomorphic to the sum of the tangent bundle 
$T{\partial M}$
of the boundary of $M$
and the normal bundle  $\nu \partial M = \partial M \times \R$
as in the diagram
\diag{%
TM_{|\partial M}  \ar[rr]^{\cong} 
\ar[dr]
& & 
T{\partial M} \oplus \nu \partial M 
\ar[dl]
\\
& \partial M & \\
}%
The sum decomposition of $TM|\partial M$ for example is obtained
by a tubular 
neighborhood $N$ of the boundary and an explicit homeomorphism 
$N\approx \partial M \times \R_+$.
The sum decomposition yields the  projection 
\begin{equation}
p_1\colon TM_{|\partial M}  \cong T{\partial M} \oplus \nu \partial M  
\to
T\partial M.
\end{equation}

Now we define a stratified manifold $X$ and the 
stratified tangent bundle $\T X$ by induction  in the number
of strata.
\begin{defi}\label{def:differentiable}
A stratified space  $X$ with finite stratification 
$X_0\subset \dots \subset X_n = X$
and attached spaces $(M_i,A_i)$ 
is termed \emph{differentiable} 
if the following conditions hold.
The spaces $M_i$ are compact differentiable manifolds with boundary
$A_i=\partial M_i$. 
Also $X_0$ is a compact differentiable manifold with 
$\partial X_0 = \emptyset$.
Moreover 
for $i=1,\dots, n$ the restriction $f^\point_i$ as in the diagram
\ndiag{%
& M_{i-1}\minus \partial M_{i-1}  \ar@{=}[d]
\\
f_i^{-1}(X_{i-1}\minus X_{i-2}) 
\ar[r]^{f^\point_i} 
\ar@{ >->}[d] 
&
X_{i-1}\minus X_{i-2} 
\ar@{ >->}[d] \\
\partial M_i 
\ar[r]^{f_i} &
X_{i-1} \\
}%
is differentiable and a 
submersion. In particular, $f_1:\partial M_1 \to X_0$
is a submersion. We assume that 
$f^{-1}_i(X_{i-1}\minus X_{i-2})$ is dense in $\partial M_i$
so that the map $f_i$ is uniquely determined by $f^\circ _i$.

We obtain inductively the tangent family $\T X_i$ in $\vec$-$\top$ 
together with inclusions
\begin{equation*}
T(X_i\minus X_{i-1}) \subset \T X_i
\end{equation*}
as follows. For $i=0$ we have $\T X_0 = TX_0$. Assume
$T(X_{i-1}\minus X_{i-2}) \subset \T X_{i-1}$ is defined.
Then we consider the diagram in $\vec$-$\top$
\diag{
&
 T( M_{i-1}\minus \partial M_{i-1})  \ar@{=}[d]
\\
T f_i^{-1}(X_{i-1}\minus X_{i-2}) 
\ar[r]^{df^\point_i} 
\ar@{ >->}[d] 
&
\ \ \ T (X_{i-1}\minus X_{i-2}) 
\ar@{ >->}[d] \\
T \partial M_i 
\ar@{.>}[r]^{df_i} &
\T X_{i-1} \\
}%
Since the inclusion in the first column of the diagram is dense there
exists at most one $\vec$-map $df_i$ for which the diagram commutes.
In addition to the differentiability in \ref{def:differentiable} we assume 
inductively that $df_i$ exists for $i\geq 2$  and that $df_i$
is surjective. 
Then $\T X_i$ is defined by the following push-out diagram,
which exists by lemma \ref{lemma:pushout}.
\ndiag{%
 & TM_i{|\partial M_i}
\ar@{ >->}[rrr]
\ar[rd]
\ar[dddd]
\ar@/_1pc/[ldd]_{p_1}
 & &  & TM_i 
\ar[dddd]
\ar[ld]
 \\
 & &  \partial M_i 
 \ar@{ >->}[r]
 \ar[dd]_{f_i}
 & M_i 
 \ar[dd]
 \\ 
T\partial M_i
\ar@/_1pc/@{.>}[rdd]_{df_i}
& &    & & \\
 & &  X_{i-1} \ar[r] & X_i \\ 
&  \T X_{i-1}
\ar[ru]
\ar[rrr]
& &  & TM\cup_{df_i \circ p_1} \T X_{i-1} = \T X_i
\ar@{.>}[ul]^{p_{i}}
\\
}%
A differentiable stratified space is termed a \emph{stratified manifold}.
\end{defi}

We call $\T X = \T X_n$  the \emph{stratified tangent bundle}
of the differentiable stratified space $X$.
As a set $\T X$ is the 
disjoint union of $T(X_i\minus X_{i-1})$,
that is of the tangent bundles $T(M_i\minus \partial M_i)$ 
for $i=0,\dots , n$.
Moreover, for every $i$, 
the pull-back of $\T X$ over 
the inclusion of the stratum $X_i \minus X_{i-1} \subset X$
is the tangent bundle of the stratum itself.
But the 
topology of $\T X$ is not the topology of the disjoint union.

\begin{lemma}
Different choices of projections $p_1$ in the definition 
of $\T X$ yield isomorphic tangent families.
\end{lemma}
\begin{proof}
We prove it by induction. For $n=0$ it is true, because 
$\T X_0 = TX_0$. So assume that $\T X_{i-1}$ does not depend
upon the choice of $p_1$.
Let $p_1$ and $p_1'$ be two projections 
$p_1,p_1': TM_i|\partial M_i \to T\partial M_i$
corresponding to two choices of the normal bundle
\begin{equation*}
\vu \partial M_i, \vu' \partial M_i  \subset T \partial M_i|\partial M_i.
\end{equation*}
The two normal bundles are equivalent and there is a homotopy 
of bundles 
$H_t: T M_i|\partial M_i$ 
covering the identity of $\partial M_i$
such that $H_1$ is the identity,
$H_t$ is the identity whenever restricted to $T\partial M_i$ 
and an isomorphism of bundles for every $t\in I$,
and $H_0$  
sends $\vu \partial M_i$ to  $\vu' \partial M_i$.
The fat homotopy $(H_t,t)\from (T M_i|\partial M_i) \times I
\to  (T M_i|\partial M_i) \times I$ defined by $(x,t) \mapsto(H_t(x),t)$
is therefore an isomorphism of $\V$-families covering the identity of 
$\partial M_i \times I$  extending the $\V$-isomorphism $F=H_0$.
Moreover, there is  a tubular neighborhood $N\approx \partial M_i \times I $
of $\partial M$ in $M_i$, where we identify $\partial M_i$ 
with $\partial M_i \times \{0\}$ in  $N$. Because 
$T M_i|N = (TM_i|\partial M_i)\times I$, the fat homotopy
$(H_t,t)$ induces a $\V$-isomorphism on $TM_i|N$ covering the identity 
of $N$, with the further property that it is the identity
on $T M_i|(\partial M_i \times \{1\})$.
Thus  
it is possible to extend $F=H_0$ 
to a $\V$-isomorphism $\tilde F$ on $TM_i$.
By construction, $p_1'F = p_1$, and therefore the following diagram
commutes and the $\V$-equivalence $G$ exists.
\ndiag{%
TM_i | \partial M_i   \ar[rrr] 
\ar[dd]_{p'_1}
& & & TM_i \ar[dddd] \\   
&  TM_i|  \partial M_i  \ar[lu]^{F} \ar[d]^{p_1} \ar[r]
& TM_i \ar[ru]^{\tilde F} \ar[dd] \\
T\partial M_i \ar[dd]^{df_i} & T\partial M_i \ar[l]^1 \ar[d]^{df_i} \\
& \T X_{i-1} \ar[r] \ar[ld]^{1} & \T X_i \ar@{.>}[dr]^{G} \\
\T X_{i-1} \ar[rrr] & & & \T X_i \\
}%
Thus the isomorphism type of 
$\T X_i$ does not depend on the choice of the projection
$p_i$, as claimed.
\end{proof}

\begin{remark}
The map $f_i$ is proper since 
$\partial M_i$ is compact and $X_{i-1}$ is Hausdorff.
Therefore
the restriction map $f_i^\point$
\begin{equation*} 
{f^\point_i}: 
f_i^{-1}(X_{i-1}\minus X_{i-2}) 
\to
X_{i-1}\minus X_{i-2} = M_{i-1}\minus \partial M_{i-1} 
\end{equation*}
is a proper surjective submersion
and hence
by Ehresmann's theorem
it is a fiber bundle. 
\end{remark}

\begin{example}\label{ex:1}
Let $M$ be a differentiable manifold with boundary; let $N=\partial M$
and 
$h=1_{\partial M}: \partial M \to \partial M=N$.
Then, the stratified space  $M\cup_h \partial M$ coincides with $M$.
The tangent  family $p:\T M \to M$
defined above is \emph{not} isomorphic to the tangent bundle $TM\to M$
(the fibers over $\partial M$ in $\T M$ are vector spaces of dimension
$\dim \partial M$ and not $\dim M$). Actually 
the stratified tangent bundle $\T M$ has the 
quotient topology under the map $\phi: TM \to \T M$ induced by the push-out.
For example let 
$M=D^n$ the closed $n$-dimensional unit disc
in $\R^n$, and $N$ its boundary $\partial M = S^{n-1}$.
$M$ is a smooth manifold and the rays through the origin give
rise to a normal bundle of $\partial M$ in $M$.
The identity map $h=1: \partial M \to S^{n-1}$ is of course 
differentiable, and the map
$dh \circ p_1$ is nothing but the projection
$TM_{|\partial M} = S^{n-1} \times \R^n \to TS^{n-1}$  
defined by  $(x,v) \mapsto (x,v - \frac{v \cdot x}{x\cdot x}x)$,
where $(x,v)$, with $x\in S^{n-1}$ and $v\in \R^n$, orthogonal to $x$, 
are the usual coordinates for the tangent bundle
of a sphere. Thus the stratified tangent bundle $\T M$ is 
defined by the quotient map $\phi: TM \to \T M$
which is the inclusion on $TM|(M\minus \partial M)$. Moreover 
two points $(x,v)$ and $(y,w)$  of 
$TM|\partial M = S^n\times \R^n$
have the same image in $\T M$ if and only if 
$v - \frac{v \cdot x}{x\cdot x}x =
w - \frac{w \cdot y}{y\cdot y}y$.
\end{example}

\begin{example}
Let $X$ be the standard cube in $\R^3$.
Let $X_0$ be the set of vertices of $X$, $X_1$ the union of the 
edges, $X_2$ the faces and $X_3=X$. Thus the filtration 
\begin{equation*}
X_0 \subset X_1 \subset X_2 \subset X_3=X
\end{equation*}
yields a decomposition into strata of dimension $0$,$1$, $2$ and 
$3$ respectively. 
It is not difficult to see that for every $i=1\dots 3$ there is a 
disjoint union of discs $M_i=\cup D^i$ and a surjective map
$f_i: \partial M_i = \cup S^{i-1} \to X_{i-1}$ such that 
$X_i$ is obtained as attachment of $M_i$ to $X_{i-1}$ via $f_i$.
The stratum $X_0$  consists of the 8 vertices,
the attaching manifold $M_1$ of 12 edges (and the boundary
of $M_1$ consists of 24  disjoint points); in a similar way,
$M_2$ is the union of 6 discs, with boundary 6 circles $S^1$,
and $M_3$ is a unique disc with boundary $S^2$.
The cube $X$ is a differentiable stratified space. 
\end{example}

\section{$\V$-function spaces}

Let $Y$ be a $\V$-family.
For $V\in \V$ we denote by $Y^V$ the space of all the $\V$-maps 
$V\to Y$, endowed with the compact-open topology. 
A sub-basis of the topology 
in $V^Y$ is given by all the sets of the form
$N_{K,U} = \{f\in Y^V \st fK \subset U \}$
with $K\subset V$ compact and $U\subset Y$ open.
For every  $\V$-family $Y$ there is a standard
projection $Y^V \to \bar Y$ sending a $\V$-map $a\from V \to Y$ 
to the point $p_Ya(0)\in \bar Y$, where $0\in V$ is the  zero
of the vector space. The pre-image of a point $b\in \bar Y$
under this projection is homeomorphic to $\hom_\V(V,Y_b)$.

Given a $\V$-map $f\from X \to Y$, let $f^V$ denote the map
$f^V\from X^V \to Y^V$ defined by composing with $f$ the
$\V$-maps $V\to X$. 
Moreover for $\varphi\from V \to W$ let $X^\varphi\from X^W \to X^V$
be the map induced by $\varphi$.
Given a $\V$-map
$f\from Z\times V \to X$, its adjoint $\hat f$ is the function
$\hat f\from Z\to X^V$ defined by $\hat f(z)(v) = f(z,v)$ for every $z\in Z$
and $v\in V$.
Lemma 
(6.10)
of \cite{qq} states that 
for every space $Z$
a $\V$-function $f\from Z\times V \to X$ is continuous if and 
only if the adjoint 
$\hat f\from Z \to X^V$ is continuous.
As a consequence,
the evaluation map
$X^V\times V \to X$ which sends $(g,v)$ to $g(v)$ is continuous
and hence 
if $f\from X\to Y$
is continuous then the induced function $f^V\from X^V \to Y^V$
is continuous, and in the same way given $\varphi\from V \to W$ 
the induced map $X^\varphi\from X^W \to X^V$ is continuous.

\begin{theo}\label{theo:coro:main}
Let $(X,D)$ be a relative $\V$-complex,
where $\V$ is a closed subcategory of $\vect$.
Let $\hh: Z_n\times D^n \to X$
be its characteristic map of $n$-cells and let  $h_n$ denote
the $n$-attaching map of $X$, i.e.~the restriction of $\hh$
to $Z_n\times S^{n-1}$.
For every $n\geq 0$ and every $V\in \V$ the diagram
\ndiag{%
(Z_n\times S^{n-1})^V \ar[r]^{\hspace{24pt} h_n^V} \ar@{ >->}[d]
& X_{n-1}^V \ar[d] \\
(Z_n\times D^n)^V \ar[r]^{\hspace{24pt} \hh^V} & X_n^V\\
}%
is a push-out (in $\top$)  with $X_{-1} = D$.
Moreover, $X^V =
 \lim_{n\to \infty} X_n^V$. 
\end{theo}
\begin{proof}
See corollary  (10.16) 
of \cite{qq}. For this we need the following lemma.
\end{proof}

A family $\mathcal{K}$ of compact sets of $V$ is termed
\emph{generating} if for every
$\V$-family $Y$
the subsets
\[
N_{K,U} = \{ f\in Y^V, f(K)\subset U\}
\]
with $K\in \mathcal{K}$
and $U$ open in $Y$ are a sub-basis for the topology of the function space 
$Y^V$.

We say that a structure category $\V$ has the 
\emph{\emph{NKC}-property} (see \cite{qq}, definition 
(6.5))
if
for every object $V\in \V$
there is a generating family of compact sets
$\mathcal{K}$
such that for every $K\in \mathcal{K}$, every
$W\in \V$ and every compact subset $C\subset W$
the subspace
\[
N_{K,C} \subset \hom_\V(V,W) = W^V
\]
is compact.

\begin{lemma}\label{lemma:NKC}
Let $\V$ be a closed subcategory of $\vect$. Then $\V$ is a NKC
category.
\end{lemma}

Since a closed subcategory of a NKC category is again
a NKC category, it suffices to consider the 
case $\V=\vect$. 
Thus lemma \ref{lemma:NKC} 
is a consequence of lemma \ref{lemma:vnkc} below.

\begin{defi}
A compact $K\subset V$ containing a linear  basis of the vector space
$V$ is termed 
\emph{spanning} compact subset.
\end{defi}

\begin{lemma}\label{lemma:compact}
Let $Y$ be a $\V$-family, 
$f\from V \to Y$ a $\V$-map, $U\subset Y$ an open  subset,
$K\subset V$ a compact  subset of $V$ 
and $f\in N_{K,U}$.
Then there is a spanning compact set $K'\subset V$ such that 
$K'\supset K$, $K'$ 
and $f\in N_{K',U} \subset N_{K,U}$.
Thus the family $\mathcal{K}$
of spanning compact subsets of $V$ is  generating.
\end{lemma}
\begin{proof}
Being $f$ continuous, $f^{-1}U$ is open in $V$,
and therefore there exists a basis $B=\{b_1,b_2,\dots, b_d\}$ 
of $V$ contained in $f^{-1}U$. Then take
$K' = K\cup B$. By definition $fB\subset U$ and hence
$f\in N_{K',U}$. Because $K'\supset K$, $N_{K',U} \subset N_{K,U}$.
\end{proof}

\begin{lemma}\label{lemma:bounded}
Let $V$ and $W$ be finitely dimensional normed vector spaces,
$B=\{b_1,b_2,\dots, b_d\}$ a basis of $V$ and $R$ 
a constant. Let $\beta\from \R^d \to V$ 
be an isomorphism sending the standard basis of $\R^d$ to the 
basis $B$ of $V$, and $|\beta^{-1}|$ the norm of its inverse.
If 
$f\from V \to
W$ is a linear map 
such that $|f(b_i)|\leq R$ for $i=1,\dots, d$, then
$|f| \leq \sqrt{d} R |\beta^{-1}|$.
\end{lemma}
\begin{proof}
The composition $f\beta\from \R^d \to W$ 
has norm
\begin{equation*}
|f\beta| = \sup_{x\in \R^d} \dfrac{|f\beta(x)|}{|x|}.
\end{equation*}
Write $x\in \R^d$ as $\sum_{j} x_j e_j$,
where $e_j$ are the elements of the standard basis of $\R^d$; 
by assumption $\beta(e_j) = b_j$ for every $j=1,\dots, d$.
Now consider the inequalities
\begin{equation*}
\dfrac{|\sum_j x_j f(b_j)|}{|x|} \leq
\dfrac{\sum_j |x_j| |f(b_j)|}{\sqrt{\sum_j x_j^2}} \leq
\sqrt{d} R.
\end{equation*}
This implies that $|f\beta| \leq \sqrt{d} R$, 
and hence $|f| \leq |f\beta \beta^{-1}| \leq \sqrt{d} R |\beta^{-1}|$,
as claimed.
\end{proof}

\begin{lemma}\label{lemma:XVclosed}
If $X$ is a trivial $\vect$-family $X=W\times Y$ over an
Hausdorff space $Y$ then 
$X^V$ is a closed subspace of $\mathrm{Map}(V,X)$.
In particular, by taking $Y=*$, we obtain that 
for every $V$, $W\in \V$, the function space 
$W^V=\vect(V,W)$ is a closed subspace of $\mathrm{Map}(V,W)$.
\end{lemma}
\begin{proof}
Let $p_X$ be the projection $W\times Y \to Y$
and 
$\{b_1,\dots, b_n\}$ a basis of the vector space $V$.
First, assume that 
$f_0 \in \mathrm{Map}(V,X)$  is a map not in $X^V$ such that 
$p_Xf(V) \neq p_Xf(0)=y_0$. Then there is $y\in Y$ 
such that $y\neq y_0$ 
and $p_Xf(V) \ni y$. Since $Y$ is Hausdorff, there are 
neighborhoods $\bar U_0$ and $\bar U_1$ in $Y$ 
such that 
$y_0\in \bar U_0$, 
$y\in \bar U_1$ and $\bar U_0\cap \bar U_1 = \emptyset$.
Let use define $U_i = W \times \bar U_i \subset X$,
$K_0 = \{0\} \subset V$ 
and $K_1 = \{v\}$,
with some $v\in (p_Xf)^{-1}y$.
We have $fK_i \subset U_i$ 
and 
\[
N_{K_0,U_0} \bigcap N_{K_1,U_1} \bigcap X^V = \emptyset,
\]
hence 
$N_{K_0,U_0} \bigcap N_{K_1,U_1}$ is an open neighborhood
of $f$ contained in $\mathrm{Map}(V,X)\minus X^V$.

On the other hand, assume that $p_Xf(V) = p_Xf(0) = y_0$.
Thus $f\from V \to p_X^{-1}y_0 = W$.
Since $f\not\in X^V$, $f$ is not linear,
there exist a vector $v\in V$ and $n$ 
scalars $\lambda_i$ in the field $\K$ such that 
\[
f( \sum_{i=1}^n \lambda_i b_i ) \neq 
\sum_{i=1}^n \lambda_i f( b_i ).
\]
Since $W\times Y$ is Hausdorff,
there are disjoint open sets $U_0$ and $U_1$ 
in $W\times Y$ such that 
$f(v) \in U_0$
and $\sum_{i=1}^n \lambda_i f( b_i ) \in U_1$.
Now, the map $W^n \to W$ with sends
$n$ vectors $v_i'$ to the linear combination
$\sum_{i=1}^n \lambda_i v'_i$ 
is continuous, 
hence there are $n$ open sets $U'_i\subset X$,
neighborhoods of $f(b_i)$,
such that if $x_i \in  U'_i$  for every $i=1,\dots,n$
and $p_X(x_i)=p_X(x_j)$ for every $i,j$
then $\sum_{i=1}^n \lambda_i x_i \subset U_1$. 
Consider the compact set $K_0=\{v\}$, $K_i=\{b_i\}\subset V$, 
$i=1,\dots,n$.
Consider the open subset $A\subset \mathbf{Map}(V,X)$ defined by
\[
A = N_{K_0,U_0} \bigcap N_{K_1,U'_1} \bigcap \dots
\bigcap N_{K_n,U'_n}.
\]
It consists of all maps $g\from V\to X$  such that 
$g(v) \subset U_0$ and $g(b_i) \subset U'_i$.
By construction $f\in A$. Moreover, 
assume that $g\in A\cap X^V$.
Then 
\[
g(v) = \sum_{i=1}^n \lambda_i g(b_i).
\]
But the left hand side of the equation belongs
to $U_0$, while the right hand side belongs 
to $U_1$, hence $g$ cannot belong to $A\cap X^V$.
Thus $A\cap X^V=\emptyset$. 
\end{proof}

\begin{lemma}\label{lemma:vnkc}
If $K$ is a spanning compact subset of $V$,
$W$ is an object of $\V=\vect$ and $C$ a compact
subset of $W$ 
then $N_{K,C}$ is compact.
\end{lemma}
\begin{proof}
As in the proof 
of lemma (10.6) 
of \cite{qq}
it is easy to see that $N_{K,C}$ is closed.
Now we prove that it is equicontinuous. 
By assumption $C$ is compact and therefore 
bounded, and by assumption $K$ contains a basis
$B=\{b_1,\dots, b_d\}$; so there exists 
a constant $R\geq 0$, $R=|C|$, such that if 
$f\in N_{K,C}$ then $|f(b_j)|\leq R$ for every $j=1,\dots, d$.
But because of lemma \ref{lemma:bounded} this implies 
that there exists a constant $R'$ (depending only
on $B$ and $R$) such that $|f| \leq R' $.
This implies  that  
for each $v_0\in V$  there is a constant $\alpha$ such that 
for every $f\in N_{K,C}$  and for every $v\in V$
\[
|f(v) - f(v_0)| \leq \alpha|v-v_0|,
\]
hence $N_{K,C}$ is equicontinuous.
Now consider for every $v_0\in V$ 
the space $\omega_{v_0} = \{ f(v_0) \}$.
By the same reason it is bounded in $W$,
hence it has compact closure in $W$.
Thus, by the Ascoli-Arzel\`a theorem,
the closure of 
$N_{K,C}$ in $\mathbf{Map}(V,W)$
is compact.
But since $N_{K,C}$ is a closed subset
of $\vect(V,W)$ which
in turn is closed
in $\mathbf{Map}(V,W)$ by lemma \ref{lemma:XVclosed},
it coincides with its closure in $\mathbf{Map}(V,W)$
and therefore is compact.
\end{proof}

With lemma \ref{lemma:vnkc} we have actually proved lemma
\ref{lemma:NKC}.

\begin{remark}
The principal bundle theorem (6.17) 
in \cite{qq}
shows that if $\V$ is a NKC category then the diagram
$X^\circ\from \V^\op \to \top$ 
given by $X^\circ(V) = X^V$
is the ``stratified principal bundle''
associated to $X$.
Moreover $X$ is determined by $X^\circ$
since $X=X^\circ \otimes_\V \ff$ 
where $\ff\from \V \to \top$ 
is the inclusion functor.
\end{remark}

\section{Algebraic constructions}

In the following we will apply theorem 
\ref{theo:coro:main}.

Let $\V$ and $\W$ be two closed subcategories of $\vect$  and let 
$\psi\from  \V \to \W$ be a continuous functor (that is,
for every $V_1,V_2\in \V$ the induced map on hom-spaces
$\psi\from  \hom_\V(V_1,V_2) \to \hom_\W(\psi V_1, \psi V_2)$ is continuous).
Given a $\V$-complex $X$ we want to build an associated 
$\W$-complex $\psi X$ over $\bar X$ with fibers
given by the images $\psi X_b$ of the fibers of $X$ under $\psi$.
At a set level, we know that  $X$ is the disjoint union 
\begin{equation*}
X = \coprod (Z_n\times e^n ),
\end{equation*}
where $e^n = D^n \minus S^{n-1}$ is the open cell
and $Z_n$ is the discrete $\V$-family of $n$-cells of $X$.
Therefore it is possible to define $\psi X$ as a set by
\begin{equation*}
\psi X = \coprod (\psi (Z_n)\times e^n ).
\end{equation*}
We now define the topology of $\psi X$ inductively.
Furthermore, we show that 
for every $V\in \V$ the map $\psi^V\from X^V \to (\psi X)^{\psi V}$
defined by sending $a\from V\to X_b\subset X$ in $X^V$ 
to $\psi(a)\in \hom(\psi(V),\psi(X_b))$, 
$\psi(a)\from \psi(V) \to \psi(X_b)=\psi(X)_b$,
is continuous.
The topology of $\psi X_0$ is the topology of   the 
$\W$-set $\psi X_0 = \psi Z_0$.
Moreover, since $\psi$ is a 
continuous functor we see that  for every $V\in \V$ the induced map
$\psi_0^V\from X_0^V \to (\psi X_0)^{\psi V}$ 
is continuous.
Now, assume that the topology of $\psi X_{n-1}$ is defined
and that
the map 
$\psi_{n-1}^V\from X_{n-1}^V \to (\psi X_{n-1})^{\psi V}$ 
is continuous.
The adjoint of the attaching map $h_n: Z_n\times S^{n-1} \to X_{n-1}$ 
is a continuous function
\begin{equation*}
\hat h_n\from S^{n-1} \to X_{n-1}^{Z_n} = \prod_{V\in Z_n} X_{n-1}^V
\end{equation*}
where the right hand side is the product space.
By inductive assumption 
$\psi_{n-1}^{Z_n}\from X_{n-1}^{Z_n} \to (\psi X_{n-1})^{\psi Z_n}$ 
is continuous because it is continuous on the factors of
$X_{n-1}^{Z_n} = \prod_{V\in Z_n} X_{n-1}^V$.  
Therefore the composition 
\ndiag{%
S^{n-1} \ar[r]^{\hat h_n} &  X_{n-1}^{Z_n} 
\ar[r]^{\psi^{Z_n}} &
 (\psi X_{n-1})^{\psi Z_n} \\
}%
is continuous and yields a $\W$-map
$\psi(h_n)\from \psi Z_n \times S^{n-1} \to \psi X_{n-1}$.
We can therefore endow $\psi X_n$ with the topology of the push-out
\ndiag{%
\psi Z_n \times S^{n-1} \ar[r]^{\hspace{20pt} \psi(h_n)} 
\ar@{ >->}[d] 
& \psi X_{n-1} \ar[d] \\
\psi Z_n \times D^n \ar[r] & \psi X_n \\
}%
Now we need to show that $\psi_{n}^{V}: X_{n}^V \to (\psi X_n)^{\psi V}$
is continuous for every $V\in \V$.
Consider the diagram \ref{diag:mydiag}.
\ldiag{diag:mydiag}{%
(\psi Z_n \times S^{n-1})^V 
\ar[ddd] \ar[rrr]^{\psi(h_n)^{\psi V}}
& & & \psi X_{n-1}^{\psi V} 
\ar[ddd]
\\
& (Z_n\times S^{n-1})^V 
\ar[ul]^{\hat \psi}
\ar[d]
\ar[r]^{\hspace{12pt} h_n^V}
& X_{n-1}^V 
\ar[ru]^{\hat \psi_{n-1}^V}
\ar[d]
& \\
& (Z_n\times D^n)^V 
\ar[r]
\ar[ld]
& 
X_{n}^V 
\ar@{.>}[rd]^{\hat \psi_n^V}
& \\
(\psi Z_n \times D^n)^V \ar[rrr] & & & \psi X_{n}^{\psi V} \\
}%
Because of theorem \ref{theo:coro:main} the two square diagrams are push-out
diagrams
and therefore the induced map $\hat \psi_n^V$ is continuous.

This inductive process yields a topology for every $\psi X_n$.
Let $\psi X = \lim_{n\geq 0} \psi X_n$. Because $\hat \psi_n^V\from
X_n^V \to (\psi X_n)^{\psi V}\subset (\psi X)^{\psi X}$  
is continuous  for every $V$, the induced map
$\hat \psi\from X^V \to (\psi X)^{\psi V}$ is continuous for 
every $V\in \V$; here we use again theorem \ref{theo:coro:main}.
Let $\VCW$ be the full subcategory of $\Vtop$
consisting of $\V$-complexes.

\begin{lemma}\label{lemma:propo:oper}
Let $\psi$ be a continuous functor $\V \to \W$ 
where $\V $ and $\W$ are closed subcategories of $\vect$.
The operation that sends each $\V$-complex $X$
to the $\W$-complex $\psi X$ defined above yields
a well-defined 
functor
$\psi\from  \VCW \to \WCW$.
\end{lemma}
This, in fact, is an example of the functor $\varphi_\#$ in 
\cite{qq} which carries a stratified bundle to an associated stratified
bundle.
\begin{proof}
We only need to define the image $\psi(f)\from \psi X \to \psi Y$ 
under $\psi$ of a $\V$-map
$f\from X\to Y$.
We define  $\psi(f)$ 
inductively on the skeleta of $\psi X$. 
Let $f_n\from X_n \to Y$ denote the restriction of 
$f$ to the $n$-skeleton of $X$.
On the zero-skeleton $\psi X_0$, $\psi (f) = \psi(f_0)$
is defined.  
Assume that $\psi (f_{n-1})$ is defined.
Let $\hh_n: Z_n \times D^n \to X_n$ be 
the $n$-characteristic map of $X$ and consider the 
adjoint of the composition 
\begin{equation*}
{(f_n\hh_n)}\from Z_n \times D^n \to X \to Y
\end{equation*}
\begin{equation*}
\widehat{(f_n\hh_n)}\from D^n \to Y^{Z_n}.
\end{equation*}
By construction the map $\hat \psi\from Y^{Z_n} \to
(\psi Y)^{\psi Z_n}$ is continuous.
Therefore the composition 
$\hat \psi \widehat{(f_n\hh_n)}$ is the adjoint
of a continuous $\W$-map
$\psi Z_n\times D^n \to \psi Y$.
Furthermore, it is possible to check that this map fits
into the outside square of the following diagram
making it commutative.
\ndiag{%
\psi Z_n \times S^{n-1} \ar[r] \ar[d] & 
\psi X_{n-1} \ar[d] \ar@/^/[rdd]^{\psi(f_{n-1})} \\
\psi Z_n \times D^n \ar[r] \ar@/_/[rrd] & 
\psi X_n \ar@{.>}[rd]_{\psi(f_n)}\\
& & \psi Y\\
}%
Hence the $\W$-map $\psi(f_n)$ exists for every $n\geq 0$.
By taking the limit, 
we obtain the desired $\W$-map $\psi(f)\from \psi X \to \psi Y$.
\end{proof}

Now consider a continuous bifunctor 
$\psi\from  \V \times \V' \to \W$, a $\V$-complex $X$ over $\bar X=B$
and a $\V'$-complex  $X'$ over the same base space $\bar X'=B$; 
assume that the cellular decompositions on $B$ induced by $X$  and $X'$
are the same and that $\V$, $\V'$ and $\W$ are closed.
We are now going to define a $\W$-complex
$\psi(X,X')$ with base space $B$.
If $X=X_0$ and $X'=X'_0$ are a $\V$-set and a $\V'$-set 
respectively over the same base set $B$,
then $\psi(X,X')$ is defined pointwise as $\psi(V,V')$,
whenever $V$ and $V'$ are fibers over the same point in $B$.
In general, let $\hh$ and $\hh'$ be the characteristic maps of $X$ and $X'$.
Using the decomposition into open cells of $X$ and $X'$ 
it is possible to define $\psi(X,X')$ as a set as follows:
\begin{equation*}
X = \coprod_{n\geq 0} Z_n \times e^n, 
\; \; X' = \coprod_{n\geq 0} Z'_n \times e^n,
\end{equation*}
\begin{equation}
\label{eq:dp}
\psi(X,X') =
\coprod_{n\geq 0} \psi(Z_n,Z_n') \times e^n. 
\end{equation}
We need to endow $\psi(X,X')$ with a suitable topology.

Consider a space $Y$ and a space $Y'$ which are 
sources of maps to the space $B$. We denote by
$Y\times_B Y'$ the pull-back
\ndiag{%
Y\times_B Y' \ar[r]\ar[d] & Y' \ar[d] \\
Y \ar[r] & B.\\
}%
For every $V\in \V$ and $V'\in \V'$, consider 
the pull-back $X^V\times_B {X'}^{V'}$ with 
respect to the maps $X^V \to \bar X$ and ${X'}^{V'} \to \bar X$. 
The pre-image of a point $b\in \bar X=\bar X'$ 
under the natural projection $X^V\times_B {X'}^{V'} \to \bar X$
is equal to $\hom(V,X_b)\times \hom(V',X'_b)$. 
Therefore we obtain a function
\begin{equation*}
\hat \psi\from 
X^V\times_B {X'}^{V'} 
\to
\psi(X,X')^{\psi(V,V')}
\end{equation*}
which is continuous on fibers (by continuity of $\psi$) defined by
\begin{equation*}
\hat \psi\from 
\hom(V,X_b)\times \hom(V',X'_b) 
\to
\hom(\psi(V,V'), \psi(X_b,X'_b)).
\end{equation*}
We want to give a topology to $\psi(X,X')$ so that $\hat \psi$
is continuous.
We can do it inductively on the skeleta of $X$ and $X'$.
If  $n=0$, then
$\psi(X_0,X_0')$ is defined as above,
and it is trivial to see that $\hat \psi_0\from X_0^V\sqcap {X_0'}^{V'} 
\to \psi(X_0,X'_0)^{\psi(V,V')}$ is continuous.
More directly, we can obtain $\psi(X,X')$ by
\[
\psi(X,X') = \psi_\#(X\times_BX')
\]
where the right hand side is defined by the associated 
bundle in \cite{qq}.

\begin{remark}\label{rem:1}
If $X$ and $X'$ are $\V$-vector bundles 
over a CW-complex $\bar X$, then 
by the bundle theorem (5.1)
of \cite{qq} $X$ and $X'$ have 
a structure of $\V$-complexes.
In this case it is not difficult to show that 
the bundles $\psi(X)$, $\psi(X,X')$ coincide 
with the corresponding constructions 
on vector bundles (see e.g. Atiyah \cite{at}, page 6).
\end{remark}

\section{Algebraic constructions on stratified vector bundles}

Consider the closure $\bar \V$ in  $\vect$ of a subcategory
$\V \subset \vect$.
Since composition maps are continuous, the closure
$\bar \V$ is a closed subcategory of $\vect$.
We say that a continuous functor $\psi\from \V \to \W$ 
is \emph{regular} if there is a continuous
functor $\bar \psi\from \bar \V \to \bar \W$
such that $\psi$ is the restriction of $\bar \psi$.
A similar definition holds for continuous bifunctors.
As a slight generalization of lemma \ref{lemma:propo:oper}
we get.

\begin{lemma}\label{propo:oper}
Let $\psi$ be a regular continuous functor $\V \to \W$ 
where $\V $ and $\W$ are subcategories of $\vect$.
Then the operation that sends each $\V$-complex $X$
to the $\W$-complex $\psi X$ defined above yields
a well-defined 
functor
$\psi\from  \VCW \to \WCW$.
\end{lemma}
\begin{proof}
The embeddings of categories $\V \subset \bar \V$
and $\W \subset \bar \W$
induce embeddings
$i_\V\from \VCW \subset  \bVCW$
and $i_\W\from \WCW \subset \bWCW $
By lemma \ref{lemma:propo:oper}, 
the continuous functor  $\bar \psi$ yields a functor
$\bar \psi\from \bVCW \to \bWCW$.
Since the image of the composition of functors $\bar \psi i_\V$
is actually in $\WCW$,
we obtain a functor $\psi\from \VCW \to \WCW$ as claimed.
\end{proof}

Let $\VCW\sqcap\VpCW$ denote the category with 
as objects pairs $(X,X')$ where $X$ is a $\V$-complex
and $X'$ a $V'$-complex, such that $\bar X$ and $\bar X'$ 
coincide. A morphism
in $\VCW\sqcap\VpCW$ 
from $(X,X')$ to $(Y,Y')$
is a pair of  maps $(f,f')$, where
$f\from X \to Y$ is a $\V$-map and $f'\from X'\to Y'$
is a $\V'$-map, such that the induced maps on the base spaces
coincide: $\bar f = \bar f'\from \bar X \to \bar Y$.
Given an object $(X,X')$ in $\VCW\sqcap\VpCW$
the $\W$-complex $\psi(X,X')$ has been defined above.
Consider a morphism $(f,f')\from (X,X') \to (Y,Y')$.
This yields a map $f\sqcap f'\from X\sqcap X' \to Y \sqcap Y'$.

More generally than \ref{eq:dp}  we get:

\begin{lemma}
\label{propo:bifunctor}
Let $\V$, $\V'$ and $\W$ be subcategories of $\vect$.
For every regular continuous bifunctor  $\psi$
$\psi\from \V\times \V' \to \W$ 
the operation that sends $(X,X')$ to $\psi(X,X')$ and 
$(f,f')\from (X,X') \to (Y,Y')$ to $\psi(f,f')$ 
is a well-defined functor 
\begin{equation*}
\psi\from 
\VCW\sqcap \VpCW \to \WCW.
\end{equation*}
\end{lemma}
\begin{proof}
The proof is identical as the proof of lemma
\ref{propo:oper}.
\end{proof}

\begin{theo}
\label{theo:construction1}
Let $\psi\from \V \to \W$ 
be a regular continuous functor, where $\V$ and $\W$ are 
subcategories
of $\vect$.
There is an induced functor $\psi\from \VStra \to \Wstra$
which sends a $\V$-stratified vector bundle $X$ to the space 
$\psi X$ over $\bar X$ with fibers $\psi X_b$,
for $b\in \bar X$.
\end{theo}
\begin{proof}
By lemma \ref{propo:oper} such functor 
is already defined for $\V$-complexes.
Using the same argument as in the proof of lemma
\ref{propo:oper}, it suffices to 
consider the case of closed (and hence NKC) subcategories
of $\vect$.
Consider the subcategories $\VStra_n$ and 
$\WStra_n$ of $\V$-stratified (resp. $\W$-stratified)
vector bundles of length at most $n$ (i.e. with  at most $n$ strata).
If $X$ is an object in $\Vstra_1$ then it is a vector bundle
over a CW-complex, hence it is a $\V$-complex by 
remark \ref{rem:1}. Thus the objects
of $\Vstra_1$ and $\WStra_1$ are in $\VCW$ and $\WCW$ 
respectively, and the maps are of course $\V$-maps 
and $\W$-maps. 
This implies that 
\[
\psi\from \Vstra_1 \to \Wstra_1
\]
is well-defined.
Assume by induction that 
$\psi$ is defined on $\Vstra_{n}$,
for some $n\geq 2$.
Let $X$ be an object in $\Vstra_{n+1}$. Thus 
there are 
CW-pairs $(\bar M_n, \bar A_n)$ 
a $\V$-vector bundle $M_n$ and a $\V$-map
$h_n\from A_n =M_n|\bar A_n \to X_{n-1}$
such that $X=X_{n-1}\cup_{h_n} M_n$.
Since $X_{n-1}$ and $A_n$ are 
both objects 
of $\Vstra_{n}$
we have that $\psi(h_n)\from \psi A_n \to \psi X_{n-1}$
is a well-defined $\W$-map.
Thus we can define $\psi X$ as the push-out 
$\psi X = \psi X_{n-1} \cup_{\psi(h_n)} M_n$. 
The same argument can be applied to a stratified 
map $f\from X \to X'$, where $X$ and $X'$ are 
objects in $\Vstra_n$.
Hence $\psi$ is well-defined on every $\Vstra_n$.
Now, since every object in $\Vstra$ belongs to 
$\Vstra_n$ for some $n$,
this implies that $\psi$ is a well-defined functor
from $\Vstra$ to $\Wstra$.
\end{proof}

Consider now a 
regular continuous bifunctor $\psi\from \V\times \V' \to \W$.
Let $\Vstra\sqcap\Vpstra$  denote the category 
of pairs $(X,X')$, where $X$ is an object in $\Vstra$ 
and $X'$ an object in $\Vpstra$ with the 
same base space $\bar X=\bar X'$.
A morphism in $\Vstra\sqcap\Vpstra$ 
is a pair of maps $(f,f')$ such that 
$f$ is a morphism in $\Vstra$ and 
$f'$ a  morphism in $\Vpstra$ which induce the same map
on the base space $\bar f = \bar f'\from \bar X \to \bar X$.

\begin{theo}\label{theo:construction2}
Every regular continuous bifunctor $\psi\from \V\times \V' \to \W$
induces a  well-defined functor 
\begin{equation*}
\psi\from 
\Vstra\sqcap \Vpstra \to \Wstra,
\end{equation*}
which sends a pair of $\V$-stratified vector bundles $X$, $X'$
over $\bar X$ to the space 
$\psi (X,X')$ over $\bar X$ with fibers $\psi (p_X^{-1}b, p_{X'}^{-1}b)$,
for $b\in \bar X$.
\end{theo}
\begin{proof}
First, assume that the structure categories are closed
in $\vect$.
As in the proof of theorem \ref{theo:construction1},
let $\Vstra_n\sqcap \Vpstra_n$ denote the subcategory
of pairs of $\V$-stratified vector bundles of length at most $n$.
The same for $\Wstra_n$.
We are going to show that 
$\psi$ is defined on $\Vstra_n\sqcap \Vpstra_n$
for every $n\geq 1$.
If $n=1$ this is true as a consequence of lemma
\ref{propo:bifunctor}.
So, assume that $\psi$ is well-defined on 
$\Vstra_{n-1}\sqcap \Vpstra_{n-1}$.
Let $(M_n,A_n)$ be the attached pair for $X$,
with attaching map $h\from A_n \to X_{n-1}$.
Let $(M'_n,A'_n)$ and $h'_n$ the corresponding for 
$X'$.
Since $\bar X = \bar X'$, 
without loss of generality we can assume that 
$\bar M_n = \bar M_n'$,
$\bar A_n = \bar A_n'$ and that 
the induced maps
$\bar h_n = \bar h_n'$ coincide.
Thus the induced map
\[
\psi(h_n,h_n')\from 
\psi(A_n,A_n') \to 
\psi(X_{n-1}, X_{n-1}')
\]
is well-defined, since $(X_{n-1},X_{n-1}')$ and $(A_n,A_m')$
belong
to $\Vstra_{n-1}\sqcap \Vpstra_{n-1}$.
The same for the inclusions
$(i,i')\from (A_n,A_n') \to (M_n,M_n')$.
Thus it is possible to define $\psi(X,X')$
as the push-out
\[
\psi(X,X') = \psi(X_{n-1},X_{n-1}') \cup_{\psi(h_n,h_n')} 
\psi(M_n,M_n').
\]
In the same way, given a pair of stratified maps
$(f,f')\from (X,X') \to (Y,Y')$ in
$\Vstra_{n}\sqcap \Vpstra_{n}$,
it is possible to define 
a $\W$-map
\[
\psi(f,f')\from 
\psi(X,X') \to \psi(Y,Y').
\]
Now, if the categories are not close,
it suffices to follow the same argument as in the proof of
lemma \ref{propo:oper}.
\end{proof}

\section{$K$-theory of stratified spaces}
Assume that the structure category $\V$
is 
closed under direct sum
$\oplus$ of vector spaces.
Then, if $\bar X$ is any space, the set $\V(\bar X)$ of isomorphism
classes of $\V$-vector bundles over $\bar X$ 
is an abelian semigroup where the sum is 
the sum induced as in lemma \ref{propo:bifunctor}
by the bifunctor $\oplus$. 
Let $K_0^\V(\bar X) = K_0(\V(\bar X))$ be the Grothendieck group
of the abelian semigroup ($\V(\bar X)$,$\oplus$).
If $\V$ is also closed under tensor product $\otimes$ 
then $K_0^\V(\bar X)$ is a ring.
This is the Atiyah-Hirzebruch $K$-theory of a space $\bar X$.

We can do the same for the category $\VStra(\bar X)$
of $\V$-stratified bundles 
over the stratified base space $\bar X$.
Let $\bar X$ be a stratified space in $\stra$. 
Then 
by theorem \ref{theo:construction2}
the category 
$\Vstra(\bar X)$
of $\V$-stratified vector bundles over $\bar X$ 
and stratified $\V$-maps over $1_{\bar X}$
is 
an additive  category with  sum $\oplus$. 

\begin{defi}\label{defi:K}
Let $\KK_0^\V(\bar X)$ denote  the Grothendieck
group $K_0(\Vstra(\bar X))$.
\end{defi}

If $\V$ is closed under the tensor product $\otimes$ then
$\KK_0^\V(\bar X)$ is a ring.
Now consider the full subcategory $\stra_0 \subset \stra$ 
of finite stratified spaces with attached spaces which are 
locally finite and countable CW-complexes.

\begin{theo}\label{theo:maintheo}
For every subcategory $\V$ of $\vect$ which is closed under 
direct sum $\oplus$ there is a functor 
\[
\KK_0^\V\from \stra^\op_0 \to \Ab, 
\]
where $\Ab$ denotes the category of Abelian groups.
If $\V$ is closed under the tensor product $\otimes$
then
the functor has its image in 
the category of rings.
Moreover there are natural homomorphisms
\ndiag{%
K_0^\V(\bar X) \ar[r] &
\KK_0^\V(\bar X) \ar[d] \ar[r] &
K_0^\V(\bar M_i) \\
& 
K_0^\V(\bar X_0)
}%
where $\bar M_i$, $\bar X_0$ are the attached spaces of $\bar X \in \Stra_0$.
\end{theo}
\begin{proof}
Let $\bar X$, $\bar X'$ be two stratified spaces in $\stra_0$ 
and $X' \to \bar X'$ an object of $\Vvect(\bar X')$.
By theorem \ref{theo:pullback}, given a stratum-preserving
map $\bar f\from \bar X \to \bar X'$ in $\Stra_0$,
the pull-back $X=\bar f^* X'$ is a well-defined 
object of $\Vvect(\bar X)$.
Moreover, a morphism $\alpha$ in $\Vvect(\bar X')$
(that is, a self-map $\alpha\from X' \to X'$ 
over the identity of $\bar X'$) induces a 
morphism $\alpha^*\from X \to X$
over the identity of $\bar X$. Hence
$\bar f^*\from 
\Vvect(\bar X') \to \Vvect(\bar X)$
is a functor.
To show that this construction yields a 
homomorphism 
$\KK_0^\V(\bar f)\from \KK_0^\V(\bar X') \to \KK_0^\V(\bar X)$
is sufficient to show that 
$\bar f^*$, as a functor, 
preserves the sum $\oplus$,
i.e. that $Y$ and $Y'$ are in $\Vvect(\bar X')$ 
then
\begin{equation}\label{eq:induction}
\bar f^*(Y\oplus Y') = \bar f^*Y \oplus
 \bar f^* Y'.
\end{equation}
Let $Y_0 \subset  Y_1 \subset \dots $ 
and $Y'_0 \subset Y_1' \subset \dots$ be 
the filtrations of $Y$ and $Y'$.
We prove by induction that  equation \ref{eq:induction}
holds for every stratified space $\bar X'$ 
with at most $n$ strata. If $n=1$ then it is true,
since in this case the $\V$-stratified bundles
are $\V$-vector bundles. So assume that 
\ref{eq:induction} is true for spaces with at 
most $n-1$ strata.
Let $(\bar M'_n,\bar A'_n)$ be 
the pair attached via the map $\bar h'$ to $\bar X'_{n-1}$
to obtain $\bar X'$,
where $\bar X'$ 
has $n$ exactly strata; in the same way,
let $(\bar M_n, \bar A_n)$ the corresponding CW-pair  
for $\bar X$ and $\bar h\from \bar A_n \to \bar X_{n-1}$
the attaching map.
Let $(M'_n,A'_n)$, $(M_n'',A''_n)$ be the 
attached bundles to $Y_{n-1}$ and $Y'_{n-1}$.
Since $\bar f$ is a map in $\stra_0$,
there is a map $\bar g\from \bar M \to \bar M'$
such that 
$\bar f=\bar f_{n-1}\cup \bar g$.
By proposition (7.4) 
of \cite{qq} we have the equality
\begin{equation}\label{eq:compare}
\bar f^*(Y\oplus Y') = 
\bar f^*(Y_{n-1}\oplus Y_{n-1}) \cup_{h} \bar g^*(M_n\oplus M'_n),
\end{equation}
where $h\from \bar g^*(A_n\oplus A'_n) \to \bar f^*(Y_{n-1}\oplus Y'_{n-1})$
is induced by the pull-back.
Now, by induction 
\[
\bar f^*(Y_{n-1}\oplus Y_{n-1}) 
=
\bar f^*Y_{n-1}\oplus \bar f^*Y_{n-1}
\]
\[
\bar g^*(M_n\oplus M'_n) =
\bar g^*M_n\oplus \bar g^*M'_n, 
\]
\[
\bar g^*(A_n\oplus A'_n) =
\bar g^*A_n\oplus \bar g^*A'_n. 
\]
But by definition the sum of the pull-backs
$\bar f^*Y\oplus \bar f^*Y'$
is equal to the push-out
\begin{equation}\label{eq:compare2}
\begin{split}
\bar f^*Y\oplus \bar f^*Y' & =
\bar f^*Y_{n-1}\oplus \bar f^*Y_{n-1} 
\cup_h
\bar g^*M_n\oplus \bar g^*M'_n \\ 
& =
\bar f^*(Y_{n-1}\oplus Y_{n-1}) \cup_{h} \bar g^*(M_n\oplus M'_n).
\end{split}
\end{equation}
Comparing equations \ref{eq:compare} and \ref{eq:compare2}
yields the desired result.
\end{proof}

\begin{theo}\label{theo:homot}
Consider a $\V$-stratified vector bundle $X'$,
a compactly stratified space $\bar X$ and a homotopy of  
(stratum-preserving) maps $\bar f_t\from \bar X \to \bar X'$,
$t\in I$.
Then 
\[
\bar f_0^* X' \cong \bar f_1^* X'.
\]
\end{theo}
\begin{proof}
We prove it by induction on the number of strata of $\bar X'$.
If $n=1$ then $X'$ is a $\V$-vector bundle
and hence the result is classical.
Assume that it is true for $\bar X'$ with at most $n-1$ strata.
Then if $\bar X$ has less than $n$ strata the claim is true.
So let $(\bar M_n,\bar A_n)$ be the CW-pair
attached to $\bar X_{n-1}$,
$(M'_n,A'_n)$ the CW-pair attached 
to $\bar X_{n-1}'$ and $\bar g_t\from (\bar M_n,\bar A_n) \to 
(\bar M_n',\bar A'_n)$ the homotopy such that 
$\bar f_t= \bar f_{n-1,t} \cup \bar g_t$, with $t\in I$.
Since $\bar X_0$ and  
$\bar M_i$ are compact for every $i\geq 1$,
we know that $\bar X_i$ is compact for every $i$.
This implies that, 
by induction, 
$\bar f^*_{0,n-1} X_{n-1}' \cong f^*_{1,n-1} X'_{n-1}$,
$\bar g_0^* A'_{n} \cong g_1^* A'_{n}$
and 
$\bar g_0^* M'_{n} \cong g_1^* M'_{n}$.
But by proposition (7.4) 
of \cite{qq} 
\[
\bar f_0^* X' = 
\bar f_{0,n-1}^* X'_{n-1} \cup \bar g_0^* M_n',
\]
\[
\bar f_1^* X' =
\bar f_{1,n-1}^* X'_{n-1} \cup \bar g_1^* M_n',
\]
and as a consequence 
$f_0^* X'\cong  f_1^* X'$.
\end{proof}

\begin{coro}\label{coro:homotopy}
If $\stra_c^\op$ denotes the full subcategory 
of $\stra^\op$ consisting 
of  compactly stratified
spaces and stratum-preserving maps,  the functor
\[
\KK_0^\V\from \stra_c^\op \to \Ab
\]
induces a functor on the homotopy category
\[
\KK_0^\V\from \stra_c^\op/{_\sim} \to \Ab.
\]
Hence if $\bar X$ is homotopy equivalent (in $\stra_c$)
to $\bar Y$, then \[
\KK_0^\V(\bar X) = \KK_0^\V(\bar Y).
\]
\end{coro}


\end{document}